\documentclass[twocolumn]{IEEEtran}
\usepackage{amsmath}
\usepackage{amssymb}
\usepackage{float}
\usepackage{multirow}
\usepackage{graphicx}

\newtheorem{thm}{Theorem}

\newtheorem{rem}{Remark}
\newtheorem{defn}{Definition}
\newtheorem{assumption}{Assumption}
\newtheorem{lem}{Lemma}

\floatstyle{ruled}
\newfloat{algorithm}{tbp}{loa}
\providecommand{\algorithmname}{Algorithm}
\floatname{algorithm}{\protect\algorithmname}
\usepackage{algorithm}
\usepackage{algorithmic}

\begin{document}

\title{Convergence of Message-Passing for Distributed Convex Optimisation with Scaled Diagonal Dominance}
\author{{\normalsize{Zhaorong Zhang$^1$, and Minyue Fu$^{1,2}$,  {\em Fellow, IEEE}}} 
%\author{...  
\thanks{$^1$School of Electrical Engineering and Computing, The University of Newcastle. University Drive, Callaghan, 2308, NSW, Australia.}
\thanks{$^2$School of Automation, Guangdong University of Technology, and Guangdong Key Laboratory of Intelligent Decision and
Cooperative Control, Guangzhou 510006, China.}
\thanks{
This work was supported by the National Natural Science Foundation of China (Grant Nos.~61633014 and U1701264). 
E-mails: zhaorong.zhang@uon.edu.au; minyue.fu@newcastle.edu.au.}
}
\maketitle
\begin{abstract}
This paper studies the convergence properties the well-known message-passing algorithm for convex optimisation. Under the assumption of pairwise separability and scaled diagonal dominance, asymptotic convergence is established and a simple bound for the convergence rate is provided for message-passing. In comparison with previous results, our results do not require the given convex program to have known convex pairwise components and that our bound for the convergence rate is tighter and simpler. When specialised to quadratic optimisation, we generalise known results by providing a very simple bound for the convergence rate. 
\end{abstract}


\section{Introduction}\label{sec1}

This paper is concerned with solving the distributed optimisation problem
\begin{align}
\min_{x\in \mathbb{R}^n} F(x) \label{eq:co}
\end{align}
for an objective function $F: \mathbb{R}^n \rightarrow \mathbb{R}$, which is assumed to be twice continuously differentiable, strictly convex and coercive (i.e., $F(x)\rightarrow \infty$ as $x\rightarrow \infty$). 

The purpose of this paper is to study the convergence properties of a well-known {\em message-passing} algorithm \cite{Pearl} for distributed optimisation, also known as {\em loopy belief propagation} \cite{Pearl}-\cite{Weiss}, {\em min-sum} \cite{Roy1}-\cite{Roy2} and {\em sum-product} \cite{Kschischang} in the literature. This paper is inspired by the excellent work of Moallemi and Van Roy~\cite{Roy1}-\cite{Roy2}, and also influenced by the great work of Malioutov, Johnson and Willsky \cite{Malioutov}, Weiss and Freeman \cite{Weiss}, and Su and Wu \cite{Su}-\cite{Su1}. 

Despite its versatile applications in many scientific and engineering disciplines, including its iconic success in error-correct decoding for approaching the Shannon coding capacity \cite{turbo}-\cite{LDPC}, the theoretic behaviour of the message-passing algorithm for loopy graphs allures many researchers for several decades. 

There has been some breakthrough recently on the convergence analysis of the message-passing algorithm for convex optimisation problems under the assumption of diagonal dominance. For quadratic optimisation, it was shown in \cite{Weiss} that the message-passing algorithm (also known as {\em Gaussian belief propagation} for marginal distribution computation of a Gaussian graphical model) convergences asymptotically to the correct minimiser (which corresponds to the correct marginal means).  This result was generalised in \cite{Malioutov} using the notion of {\em walk-summability} to show that the same asymptotic convergence is guaranteed  under a relaxed assumption of {\em generalised diagonal dominance} (or {\em scaled diagonal dominance \cite{Roy2})}. In \cite{Roy2}, the same convergence property for quadratic optimisation is shown under an equivalent assumption of {\em pairwise convex separability} and more flexible initial messages.  General necessary and sufficient conditions for asymptotic convergence of the message-passing algorithm for quadratic optimisation were established in \cite{Su}-\cite{Su1}, but verification of these conditions can be difficult. In \cite{Roy1}, convergence properties of  message-passing were generalised to convex optimisation. For a {\em pairwise separable convex program} with scaled diagonal dominance, asymptotic convergence and a bound for the convergence rate were established.

This paper is motivated by the fact that the convergence analysis in \cite{Roy1} requires a strong pairwise separation form, i.e.,  {\em every} pairwise component must be {\em convex}. This assumption is not consistent with the results for quadratic optimisation in \cite{Roy2} and \cite{Malioutov} where no such constraint is required. Although it is true for quadratic optimisation that scaled diagonal dominance is equivalent to pairwise separability with convex components and that this equivalence may be extendable to the non-quadratic case, searching for such pairwise separation may constitute a separate optimisation problem.  In this paper, we follow the same analysis method as in \cite{Roy1}. Through more careful analysis, we discover that the requirement for convex components can indeed be dropped, leading to a true generalisation of convergence results for quadratic optimisation. More specifically, we show that for any strict convex program in a pairwise separation form, asymptotic convergence is guaranteed under scaled diagonal dominance. Moreover, by choosing  appropriate initial messages, we provide a bound for the convergence rate tighter and simpler than that in \cite{Roy1}. When specialised to  quadratic optimisation, our results generalise the work of  \cite{Roy2} and \cite{Malioutov} by having a very simple bound for the convergence rate. 

The rest of the paper is organised as follows. Section~\ref{sec2} formulates the distributed convex optimisation problem and introduces the message-passing algorithm; Section~\ref{sec3} presents a well-poseness result on the algorithm and our main results on convergence analysis; Section~\ref{sec4} presents the proofs; Section~\ref{sec5} specialises our results to quadratic optimisation; and Section~\ref{sec6} concludes the paper. 

\section{Problem Formulation}\label{sec2}

\subsection{Distributed Convex Optimisation}

Consider the objective function $F(\cdot)$ in (\ref{eq:co}) and denote its optimal solution by $x^{\star}$. 

\begin{defn} \label{def:ps}
({\em Pairwise Separation}) An objective function $F: \mathbb{R}^n \rightarrow \mathbb{R}$ is said to be in a {\em pairwise separation} form if it is expressed as 
\begin{align}
F(x) &= \sum_{i\in V} f_i(x_i) + \sum_{(i,j)\in E} f_{ij}(x_i, x_j), \label{eq:ps}
\end{align}
for some graph $G=(V, E)$ with node set $V=\{1, 2, \ldots, n\}$ and undirected edge set $E\subset V\times V$, and that the factors $\{f_i(\cdot)\}$ and $\{f_{ij}(\cdot, \cdot)\}$ are twice continuously differentiable and coercive.  $F(x)$ is said to be in a {\em pairwise convex separation} form if, in addition, $\{f_i(\cdot)\}$ are strictly convex and $\{f_{ij}(\cdot, \cdot)\}$ are convex. 
\end{defn}

\begin{defn} \label{def:dd}
({\em Scaled Diagonal Dominance} \cite{Roy1}) Given a twice continuously differentiable function $F: \mathbb{R}^n \rightarrow \mathbb{R}$, a scalar $\lambda\in (0, 1)$ and a positive vector $w=(w_1, w_2, \ldots, w_n)\in \mathbb{R}^n$, $F$ (or $\frac{\partial^2}{\partial x^2}F(x)$) is said to be $(\lambda, w)$-scaled diagonally dominant if, for every $i\in V$ and all $x\in \mathbb{R}^n$, 
\begin{align}
\sum_{j\in N_i} w_j\left |\frac{\partial^2}{\partial x_i\partial x_j}F(x)\right | &\le \lambda w_i  \frac{\partial^2}{\partial x_i^2}F(x). 
\label{eq:dd}
\end{align}
\end{defn}

\begin{rem}\label{rem:1}
It is obvious that pairwise convex separation is a stronger condition than pairwise separation.  It is pointed out in \cite{Roy1} that a quadratic $F(x)$ is {\em pairwise convex separable} if and only if it is scaled diagonally dominant. However, finding such pairwise convex separation may constitute a separate optimisation problem. 
\end{rem}

\begin{assumption}\label{as:1}
The given objective function $F: \mathbb{R}^n \rightarrow \mathbb{R}$ has the following properties:
\begin{itemize}
\item $F(\cdot)$ is strictly convex and coercive;
\item $F(\cdot)$ is expressed in a pairwise separation form (\ref{eq:ps}) with some graph $G=(V,E)$; 
\item $F(\cdot)$ is scaled diagonally dominant for some $0<\lambda<1$ and positive $w \in \mathbb{R}^n$.
\end{itemize}
\end{assumption}
 
\subsection{Message-Passing Algorithm}

The {\em message passing} algorithm is a distributed iterative algorithm operating on each node $i\in V$. At iteration $t=0, 1, \ldots$, each node $i$ takes an {\em incoming message} $J_{u\rightarrow i}^{(t)}: \mathbb{R}\rightarrow \mathbb{R}$ from each neighbouring node $u\in N_i$. These incoming messages are fused together to create an {\em outgoing message} $J_{i\rightarrow j}^{(t+1)}: \mathbb{R}\rightarrow \mathbb{R}$ for each $j\in N_i$ according to
\begin{align}
J_{i\rightarrow j}^{(t+1)}(x_j) :=& \min_{y_i} f_i(y_i) + f_{ji}(x_j,y_i) + \sum_{u\in N_i \backslash j} J_{u\rightarrow i}^{(t)}(y_i)\nonumber \\
&  + \kappa_{i\rightarrow j}^{(t+1)}, \label{eq:mp}
\end{align}
where $\kappa_{i\rightarrow j}^{(t+1)}$ is such that $J_{i\rightarrow j}^{(t+1)}(0)=0$. The local objective function for each $x_i$  in iteration $t$ is constructed as
\begin{align}
\gamma_i^{(t)}(x_i) &:= f_i(x_i) + \sum_{u\in N_i} J_{u\rightarrow i}^{(t)}(x_i), \label{eq:git}
\end{align}
and the estimate $x_i^{(t+1)}$ for $x_i^{\star}$ (the $i$-th component of $x^{\star}$) in iteration $t$ is obtained by minimising $\gamma_i^{(t)}$, i.e.,  
\begin{align}
x_i^{(t+1)}&:=\mathrm{argmin\ }  \gamma_i^{(t)}(x_i). \label{eq:xit}
\end{align}
 
We will consider the choice of initial messages as follows. 

\begin{assumption}\label{as:2}
Take any set of initial estimates $\{x_{i\rightarrow j}^{(0)}, i\in V, j\in N_i\}$ and any set of initial messages $\{J_{i\rightarrow j}^{(0)}(\cdot), i\in V, j\in N_i\}$ which are twice continuously differentiable and satisfy the following for some $0\le \rho< \lambda^{-1}$:
\begin{align}
&\frac{d^2}{dx_j^2}J_{i\rightarrow j}^{(0)}(x_j) - \frac{d^2}{dx_j^2}f_{ji}(x_j, x_{i\rightarrow j}^{(0)})\nonumber \\
\ge& -\rho\frac{w_i}{w_j}\left |\frac{\partial^2}{\partial x_j\partial x_i} f_{ji}(x_j,x_{i\rightarrow j}^{(0)})\right |,  \ \forall x_j\in \mathbb{R}.\label{eq:Ji1}
\end{align}
\end{assumption}

\begin{rem}\label{rem:3}
Note, in particular, that any  
\begin{align*}
J_{i\rightarrow j}^{(0)}(x_j)&=f_{ji}(x_j,x_{i\rightarrow j}^{(0)})+c_{ji}(x_j,x_{i\rightarrow j}^{(0)})
\end{align*}
with $c_{ji}(\cdot,\cdot)$ being affine in the first variable will satisfy (\ref{eq:Ji1}), including $c_{ji}(\cdot,\cdot)=0$. 
Also note that our requirement for the initial messages is weaker than that in \cite{Roy1} because our $J_{i\rightarrow j}^{(0)}(\cdot)$  are not necessarily convex, whereas they must be in \cite{Roy1} (due to the components $f_{ji}(\cdot,\cdot)$ being required to be convex).  
\end{rem}
 
{\bf Notation}: The symbol $\nabla$ is used to denote partial derivative, i.e., $\nabla f(x) = \frac{\partial}{\partial x}f(x)$, $\nabla^2 f(x) = \nabla (\nabla f(x))$. Similarly, $\nabla_1 f(x,y)=\frac{\partial}{\partial x}f(x,y)$, $\nabla_2 f(x,y)=\frac{\partial}{\partial y}f(x,y)$, $\nabla_{12} f(x,y)=\frac{\partial^2}{\partial x\partial y}f(x,y)$, $\nabla_1^2 f(x,y)=\frac{\partial^2}{\partial x^2}f(x,y)$, etc.. 
 
\section{Main Results}\label{sec3}

\subsection{Well-Posedness of Message-Passing}

Our first task is to determine whether the message-passing algorithm is well posed or not. By this, we mean whether or not the minimisation problems (\ref{eq:mp}) and (\ref{eq:xit}) have unique solutions, i.e., $J_{i\rightarrow j}^{(t+1)}(\cdot)$ and $x_i^{(t+1)}$ are well defined. We show below that this property is guaranteed under Assumptions~\ref{as:1}-\ref{as:2}.

For any $i\in V, j\in N_i$ and $t\ge0$, define $g_{ij}^{(t)}:\mathbb{R}\times\mathbb{R}\rightarrow \mathbb{R}$ as
\begin{align}
g_{ij}^{(t)}(x_i, x_j) &:= f_i(x_i)+f_{ji}(x_j,x_i)+\sum_{u\in N_i\backslash j} J_{u\rightarrow i}^{(t)}(x_i) \label{eq:gij}
\end{align}
and denote its minimiser by  
\begin{align}
x_{i\rightarrow j}^{(t+1)}(x_j) &= \mathrm{arg}\min_{y_i} g_{ij}^{(t)}(y_i,x_j) \label{eq:yi}
\end{align}
(or $x_{i\rightarrow j}^{(t+1)}$ for short). Also define $a_{i\rightarrow j}^{(t)}: \mathbb{R}\rightarrow \mathbb{R}$ as
 \begin{align}
 a_{i\rightarrow j}^{(t)}(x_j) := &\nabla_1^2 g_{ij}^{(t)}(x_{i\rightarrow j}^{(t+1)},x_j). \label{eq:aij}
\end{align}

\begin{thm}\label{thm:1}
Under Assumptions~\ref{as:1}-\ref{as:2}, the following properties hold for all $i\in V, j\in N_i$ and $t=0,1,\ldots$.
\begin{itemize}
\item[P1:] $\nabla^2 \gamma_i^{(t)}(x_i)>0$ and $\nabla_1^2 g_{ij}^{(t)}(x_i,x_j)>0$ for all $x_i, x_j\in \mathbb{R}$, i.e., $\gamma_i^{(t)}(\cdot)$ is strictly convex and $g_{ij}^{(t)}(\cdot,\cdot)$ is strictly convex with respect to the first variable, implying that the message-passing algorithm is well posed;
\item[P2:]  $a_{i\rightarrow j}^{(t)}(x_j)>0$ for all $x_j\in \mathbb{R}$; 
\item[P3:] $\lambda w_i a_{i\rightarrow j}^{(t)}(x_j)\ge w_j |\nabla_{12} f_{ji}(x_j,x_{i\rightarrow j}^{(t+1)})|$ for all $x_j\in \mathbb{R}$;
\item[P4:] For all $x_j\in \mathbb{R}$, it holds that
\begin{align}
\hspace{-5mm}\nabla x_{i\rightarrow j}^{(t+1)}(x_j)=&-\frac{\nabla_{12} f_{ji}(x_j, x_{i\rightarrow j}^{(t+1)})}{a_{i\rightarrow j}^{(t)}(x_j)}, \label{eq:dxij} \\
\hspace{-5mm}\nabla^2 \hspace{-1mm}J_{i\rightarrow j}^{(t+1)}(x_j)=&\nabla_1^2 f_{ji}(x_j,x_{i\rightarrow j}^{(t+1)})-\frac{|\nabla_{12} f_{ji}(x_j, x_{i\rightarrow j}^{(t+1)})|^2}{a_{i\rightarrow j}^{(t)}(x_j)}, \label{eq:Jxij} 
\end{align}
and in particular, 
\begin{align}
0&\ge \nabla^2 J_{i\rightarrow j}^{(t+1)}(x_j)-\nabla_1^2 f_{ji}(x_j,x_{i\rightarrow j}^{(t+1)})\nonumber \\
&\ge -\lambda \frac{w_i}{w_j}|\nabla_{12} f_{ji}(x_j, x_{i\rightarrow j}^{(t+1)})|. \label{eq:Jxij1} 
\end{align}
\end{itemize} 
\end{thm}

\begin{IEEEproof}
Define 
\begin{align*}
L_{ij}^{(t)}(x_i,x_j) :=& \nabla_1 g_{ij}^{(t)}(x_i,x_j) \\
=& \nabla f_i(x_i) + \nabla_2 f_{ji}(x_j,x_i) + \hspace{-2mm}\sum_{u\in N_i\backslash j} \nabla J_{u\rightarrow i}^{(t)}(x_i). 
\end{align*}
We first verify P1-P5 for $t=0$. 

Using (\ref{eq:Ji1}), we get 
\begin{align*}
&\nabla_1^2 g_{ij}^{(0)}(x_i,x_j) \\
=&\nabla_1 L_{ij}^{(0)}(x_i,x_j) \\
=& \nabla^2 f_i(x_i) + \nabla_2^2 f_{ji}(x_j,x_i) + \hspace{-2mm}\sum_{u\in N_i\backslash j} \nabla^2 J_{u\rightarrow i}^{(0)}(x_i) \\
\ge& \nabla^2 f_i(x_i) + \nabla_2^2 f_{ji}(x_j,x_i) + \hspace{-2mm}\sum_{u\in N_i\backslash j} \nabla_1^2 f_{iu}(x_i,x_{u\rightarrow i}^{(0)})\\
&-\frac{\rho}{w_i} \sum_{u\in N_i\backslash j}w_u |\nabla_{12}f_{iu}(x_i,x_{u\rightarrow i}^{(0)})|.
\end{align*}
The first three terms on the right hand of the inequality above equals to $\frac{\partial^2}{\partial x_i^2}F(x)$ with $x$ evaluated with $x_u=x_{u\rightarrow i}^{(0)}$ for all $u\in N_i\backslash j$. So,  
\begin{align*}
&\nabla_1^2 g_{ij}^{(0)}(x_i,x_j) \\
\ge&\frac{\partial^2}{\partial x_i^2} F(x) +\rho\frac{w_j}{w_i}|\nabla_{12}f_{ij}(x_i,x_j)|-\rho\frac{w_j}{w_i}|\nabla_{12}f_{ij}(x_i,x_j)|\\
&- \frac{\rho}{w_i} \sum_{u\in N_i\backslash j}w_u |\nabla_{12}f_{iu}(x_i,x_{u\rightarrow i}^{(0)})|\\
\ge& \frac{\partial^2}{\partial x_i^2} F(x)+ \rho\frac{w_j}{w_i}|\nabla_{12}f_{ij}(x_i,x_j)|-\rho\lambda \frac{\partial^2}{\partial x_i^2} F(x) \\
\ge& (1-\rho\lambda)\frac{\partial^2}{\partial x_i^2} F(x) + \rho\frac{w_j}{w_i}|\nabla_{12}f_{ij}(x_i,x_j)|>0.
\end{align*}
The second inequality above used the scaled diagonal dominance property.
This verifies the property of $g_{ij}^{(0)}(\cdot,\cdot)$ in P1. 

Subsequently, $x_{i\rightarrow j}^{(1)}(x_j)$ is well defined, and (\ref{eq:aij}) is well defined too. It follows from above inequality that P2 holds for $t=0$. Taking the above inequality further, we get
\begin{align*}
&\nabla_1^2 g_{ij}^{(0)}(x_i,x_j) \\
\ge&(1-\rho\lambda)\frac{\partial^2}{\partial x_i^2} F(x) + \rho\frac{w_j}{w_i}|\nabla_{12}f_{ij}(x_i,x_j)|\\
\ge&(1-\rho\lambda)\lambda^{-1} \frac{w_j}{w_i}|\nabla_{12}f_{ij}(x_i,x_j)|  + \rho\frac{w_j}{w_i}|\nabla_{12}f_{ij}(x_i,x_j)|\\
=&\frac{w_j}{\lambda w_i}|\nabla_{12}f_{ij}(x_i,x_j)|=\frac{w_j}{\lambda w_i}|\nabla_{12}f_{ji}(x_j,x_i)|.
\end{align*}
Taking $x_i=x_{i\rightarrow j}^{(1)}$ and using (\ref{eq:aij}), P3 holds for $t=0$. 

Next, note that $x_{i\rightarrow j}^{(1)}(x_j)$ solves $L_{ij}^{(0)}(x_{i\rightarrow j}^{(1)},x_j)=0$. It follows that 
\begin{align*}
0=&\frac{d}{dx_j}L_{ij}^{(0)}(x_{i\rightarrow j}^{(1)},x_j)\\
=&\nabla_1L_{ij}^{(0)}(x_{i\rightarrow j}^{(1)},x_j)\nabla x_{i\rightarrow j}^{(1)} + \nabla_2L_{ij}^{(0)}(x_{i\rightarrow j}^{(1)},x_j)\\
=&a_{i\rightarrow j}^{(0)}(x_j)\nabla x_{i\rightarrow j}^{(1)}+\nabla_{12} f_{ji}(x_j,x_{i\rightarrow j}^{(1)}),
\end{align*}
which gives (\ref{eq:dxij}) for $t=0$. 

Next, note that $J_{i\rightarrow j}^{(1)}(x_j) = g_{ij}^{(0)}(x_{i\rightarrow j}^{(1)}, x_j)$. It follows that
\begin{align*}
\nabla J_{i\rightarrow j}^{(1)}(x_j) &= \nabla_1 g_{ij}^{(0)}(x_{i\rightarrow j}^{(1)},x_j)\nabla x_{i\rightarrow j}^{(1)} + \nabla_2 g_{ij}^{(0)}(x_{i\rightarrow j}^{(1)},x_j) \\
&=L_{ij}^{(0)}(x_{i\rightarrow j}^{(1)},x_j)\nabla x_{i\rightarrow j}^{(1)}+\nabla_1 f_{ji}(x_j,x_{i\rightarrow j}^{(1)}) \\
&=\nabla_1 f_{ji}(x_j,x_{i\rightarrow j}^{(1)}).
\end{align*}
Differentiating the above again and using (\ref{eq:dxij}), we have
\begin{align*}
\nabla^2 J_{i\rightarrow j}^{(1)}(x_j) &=  \nabla_1^2 f_{ij}(x_j,x_{i\rightarrow j}^{(1)})+\nabla_{12} f_{ji}(x_j,x_{i\rightarrow j}^{(1)})\nabla x_{i\rightarrow j}^{(1)}\\
&=\nabla_1^2 f_{ji}(x_j,x_{i\rightarrow j}^{(1)}) -\frac{|\nabla_{12} f_{ji}(x_j,x_{i\rightarrow j}^{(1)})|^2}{a_{i\rightarrow j}^{(0)}(x_j)}.
\end{align*}
This proves (\ref{eq:Jxij}) and the first inequality in (\ref{eq:Jxij1}) for $t=0$.  Using P3 in the above equation, we further get 
\begin{align*}
\nabla^2 J_{i\rightarrow j}^{(1)}(x_j)&\ge \nabla_1^2 f_{ji}(x_j,x_{i\rightarrow j}^{(1)}) -\lambda\frac{w_i}{w_j}|\nabla_{12} f_{ij}(x_j,x_{i\rightarrow j}^{(1)})|,
\end{align*}
which is the second inequality in (\ref{eq:Jxij1}) for $t=0$. Hence, P4 holds for $t=0$. 

Next,  using Assumption~\ref{as:2} again, we have
\begin{align*}
\nabla^2 \gamma_i^{(0)}(x_i) =& \nabla^2 f_i(x_i) + \sum_{u\in N_i}\nabla^2 J_{u\rightarrow i}^{(0)}(x_i)\\
\ge& \nabla^2 f_i(x_i) + \sum_{u\in N_i} \nabla_1^2 f_{iu}(x_i,x_{u\rightarrow i}^{(0)})\\
&-\rho\sum_{u\in N_i}\frac{w_u}{w_i}|\nabla_{12}f_{iu}(x_i,x_{u\rightarrow i}^{(0)})|\\
=&\frac{\partial^2}{\partial x_i^2} F(x) - \rho\sum_{u\in N_i}\frac{w_u}{w_i}|\nabla_{12}f_{iu}(x_i,x_{u\rightarrow i}^{(0)})|\\
\ge&   \frac{\partial^2}{\partial x_i^2} F(x) - \rho\lambda \frac{\partial^2}{\partial x_i^2} F(x)>0.
\end{align*}
In the above, $x$ is evaluated with $x_u=x_{u\rightarrow i}^{(0)}$ for all $u\in N_i\backslash j$. This verifies the property of $\gamma_i^{(0)}(\cdot)$ in P1 for $t=0$.

By now, we have verified all the P1-P5 for $t=0$. The proof for $t>1$ is done by repeating the above steps. Note that (\ref{eq:Ji1}) is replaced with (\ref{eq:Jxij1}) in subsequent iterations in which $\rho$ becomes $\lambda$, but the condition $\rho<\lambda^{-1}$ remains valid. 
\end{IEEEproof}

\subsection{Convergence Properties of Message-Passing}

We now present the main results of this paper. Its proof will be presented later. 

\begin{thm}\label{thm:2}
Under Assumptions~\ref{as:1}-\ref{as:2}, the following convergence property holds for every node $r\in V$ and $t>0$:
\begin{align}
&\frac{|x_r^{(t)}-x_r^{\star}|}{w_r}\nonumber \\
\le& \frac{\lambda^{t}}{1-\lambda}\max_{i\in V} \frac{\sum_{u\in N_i}|\nabla_1 f_{iu}(x_i^{\star},x_u^{\star})-\nabla J_{u\rightarrow i}^{(0)}(x_i^{\star})|}{w_i\min_{x\in \mathbb{R}^n} \frac{\partial^2}{\partial x_i^2}F(x)}.
\label{eq:rate}
\end{align}
Moreover, if the initial messages are chosen to be $J_{i\rightarrow j}^{(0)}(x_j)=f_{ji}(x_j,x_i^{(0)})$ and $x_{i\rightarrow j}^{(0)}=x_i^{(0)}$ for all $i\in V, j\in N_i$, then (\ref{eq:rate}) simplifies to 
\begin{align}
\frac{|x_r^{(t)}-x_r^{\star}|}{w_r}&\le \lambda^{t+1}\frac{M}{1-\lambda}
\max_{v\in V}\frac{|x_{v}^{(0)}-x_v^{\star}|}{w_v}, \label{eq:rate1}
\end{align}
where $M$ is the {\em conditioning value} for the diagonal part of the Jacobian matrix for $F(\cdot)$, i.e., 
\begin{align}
M=\max_{i\in V} \frac{\max_{x\in \mathbb{R}^n} \frac{\partial^2}{\partial x_i^2}F(x)}{\min_{x\in \mathbb{R}^n} \frac{\partial^2}{\partial x_i^2}F(x)}. \label{eq:M}
\end{align}
\end{thm}

\section{Proof of Theorem~\ref{thm:2}}\label{sec4}

The basic idea of the proof follows from \cite{Roy1}. There are two main differences though: 1) We need to handle non-convex initial messages; 2) The converge bound in Theorem~\ref{thm:2} gives the direct link between the initial estimation errors and those in each iteration. Extra efforts are needed in the proof due to these differences. Similar to \cite{Roy1}, the proof relies on two critical tools: One is the use of {\em parameterised initial messages} \cite{Roy1}; Another is the {\em computation tree} \cite{Roy1} (also known as {\em unwrapped tree} in \cite{Weiss,Malioutov}.  We first introduce these tools. 
 
\subsection{Parameterised Initial Messages}

Similar to Lemma 2 of \cite{Roy1},  consider the following parameterised initial messages 
\begin{align}
J_{i\rightarrow j}^{(0)}(x_j,p) &:= J_{i\rightarrow j}^{(0)}(x_j) + x_jpC_{ji} \label{eq:initial1}
\end{align}
where the parameter $p\in [0, 1]$ and  
\begin{align}
C_{ji}:=& \nabla_1 f_{ji}(x_j^{\star},x_i^{\star})-\nabla J_{i\rightarrow j}^{(0)}(x_j^{\star}), \label{Cji}
\end{align}  
The evolution of messages  is revised to be 
\begin{align*}
J_{i\rightarrow j}^{(t+1)}(x_j,p) &=\min_{y_i} g_{ij}^{(t)}(y_i, x_j,p)+\kappa_{i\rightarrow j}^{(t+1)}
\end{align*}
with 
\begin{align*}
g_{ij}^{(t)}(x_i, x_j,p) &:= f_i(x_i)+f_{ji}(x_j,x_i)+\sum_{u\in N_i\backslash j} J_{u\rightarrow i}^{(t)}(x_i,p).
\end{align*}
Similarly, the revised local objective functions and their optimal estimates are given by 
\begin{align*}
x_i^{(t+1)}(p)&=\mathrm{arg}\min_{x_i} \gamma_i^{(t)}(x_i,p):= f_i(x_i) + \sum_{u\in N_i} J_{u\rightarrow i}^{(t)}(x_i,p).
\end{align*}

\begin{lem}\label{lem:2}
Suppose Assumptions~\ref{as:1}-\ref{as:2} hold and the  initial messages (\ref{eq:initial1}) are used. Then, for $p=1$, 
we have
\begin{align}
\nabla_1 J_{i\rightarrow j}^{(t)}(x_j^{\star},1)&=\nabla_1 f_{ji}(x_j^{\star},x_i^{\star})
\end{align}
 and $x_i^{(t+1)}(1)=x_i^{\star}$ for all $i\in V$ and all $t\ge0$.
\end{lem}

\begin{IEEEproof}
It is clear that (\ref{eq:initial1}) still satisfies Assumption~\ref{as:2} for any $p$.  We first verify the result for $t=0$. Indeed,  
\begin{align*}
\nabla_1 J_{i\rightarrow j}^{(0)}(x_j^{\star}, 1) &= \nabla J_{i\rightarrow j}^{(0)}(x_j^{\star}) + C_{ji} = \nabla_1 f_{ji}(x_j^{\star},x_i^{\star}),
\end{align*}
which leads to 
\begin{align*}
\nabla_1 \gamma_i^{(0)}(x_i^{\star},1) &= \nabla_1 f_i(x_i^{\star})+\sum_{u\in N_i} \nabla_1 J_{u\rightarrow i}^{(0)}(x_i^{\star},1) \\
&= \nabla_1 f_i(x_i^{\star})+\sum_{u\in N_i} \nabla_1 f_{iu}(x_i^{\star},x_u^{\star}) =0.
\end{align*}
This coincides with the first-order optimality condition for $x_i^{\star}$. It follows that $x_i^{(1)}(1) = x_i^{\star}$.  

Now consider $t=1$. Minimising $g_{ij}^{(0)}(y_i, x_j^{\star},1)$ yields $y_i=x_i^{\star}$ because
\begin{align*}
&\nabla_1 g_{ij}^{(0)}(x_i^{\star}, x_j^{\star},1)\\
=&\nabla f_i(x_i^{\star})+\nabla f_{ji}(x_j^{\star},x_i^{\star})+\sum_{u\in N_i\backslash j} \nabla J_{u\rightarrow i}^{(0)}(x_i^{\star},1)\\
=&\nabla f_i(x_i^{\star}) + \nabla_2 f_{ji}(x_j^{\star}, x_i^{\star})+\sum_{u\in N_i\backslash j} \nabla_1f_{ji}(x_j^{\star},x_i^{\star}) =0,
\end{align*}
which coincides with the first-order optimality condition for $x_i^{\star}$. It follows that
\begin{align*}
\nabla_1 J_{i\rightarrow j}^{(1)}(x_j^{\star},1)&= \frac{d}{dx_j^{\star}}g_{ij}^{(0)}(x_i^{\star},x_j^{\star},1) \\
&=\nabla_1 g_{ij}^{(0)}(x_i^{\star}, x_j^{\star},1)\nabla x_i^{\star} + \nabla_2 g_{ij}^{(0)}(x_i^{\star}, x_j^{\star},1)\\
&=\nabla_2 g_{ij}^{(0)}(x_i^{\star}, x_j^{\star},1)\\
&=\nabla_1 f_{ji}(x_j^{\star}, x_i^{\star}). 
\end{align*}
Like the case of $t=0$, the above yields $x_i^{(1)}(1) = x_i^{\star}$. The proof for $t>1$ is done by repeating the above process.  
\end{IEEEproof}

\subsection{Computation Tree}

\begin{figure}[ht]
\begin{picture}(240,135)
\put(54,124){\circle{6}}
\put(57,123){\line(1,-1){39}} 
\put(51,123){\line(-1,-1){39}} 
\put(54,121){\line(0,-1){38}} 
\put(54,80){\circle{6}}
\put(12,81){\circle{6}}
\put(96,81){\circle{6}}
\put(54,77){\line(0,-1){38}}
\put(56,38){\line(1,1){40}} 
\put(52,38){\line(-1,1){40}} 
\put(54,36){\circle{6}}
\put(58,30){$5$}
\put(58,126){$1$}
\put(58,75){$3$}
\put(17,75){$2$}
\put(99,75){$4$}

\put(176,134){\circle{6}}
\put(179,133){\line(1,-1){39}} 
\put(173,133){\line(-1,-1){39}} 
\put(176,131){\line(0,-1){38}} 
\put(176,90){\circle{6}}
\put(134,91){\circle{6}}
\put(218,91){\circle{6}}
\put(176,87){\line(0,-1){25}}
\put(134,88){\line(0,-1){25}}
\put(218,88){\line(0,-1){25}}
\put(134,60){\circle{6}}
\put(176,60){\circle{6}}
\put(218,60){\circle{6}}
\put(180,136){$1$}
\put(180,85){$3$}
\put(139,85){$2$}
\put(221,85){$4$}
\put(139,60){$5$}
\put(180,60){$5'$}
\put(222,60){$5''$}
\put(133,57){\line(-3,-5){12}}
\put(135,57){\line(3,-5){12}}
\put(175,57){\line(-3,-5){12}}
\put(177,57){\line(3,-5){12}}
\put(217,57){\line(-3,-5){12}}
\put(219,57){\line(3,-5){12}}
\put(121,34){\circle{6}}
\put(147,34){\circle{6}}
\put(163,34){\circle{6}}
\put(190,34){\circle{6}}
\put(205,34){\circle{6}}
\put(231,34){\circle{6}}
\put(121,31){\line(0,-1){24}}
\put(147,31){\line(0,-1){24}}
\put(163,31){\line(0,-1){24}}
\put(190,31){\line(0,-1){24}}
\put(205,31){\line(0,-1){24}}
\put(231,31){\line(0,-1){24}}
\put(121,4){\circle{6}}
\put(147,4){\circle{6}}
\put(163,4){\circle{6}}
\put(190,4){\circle{6}}
\put(205,4){\circle{6}}
\put(231,4){\circle{6}}
\put(110,32){$4'$}
\put(137,32){$3'$}
\put(167,32){$2'$}
\put(177,32){$4''$}
\put(209,32){$3''$}
\put(235,32){$2''$}
\put(110,2){$1'$}
\put(136,2){$1''$}
\put(166,2){$1'''$}
\put(179,2){$1`$}
\put(209,2){$1``$}
\put(235,2){$1```$}
\end{picture}
  \caption{Loopy graph on the left. The computation tree for node 1 on the right.}\label{fig:2-1}
\end{figure}
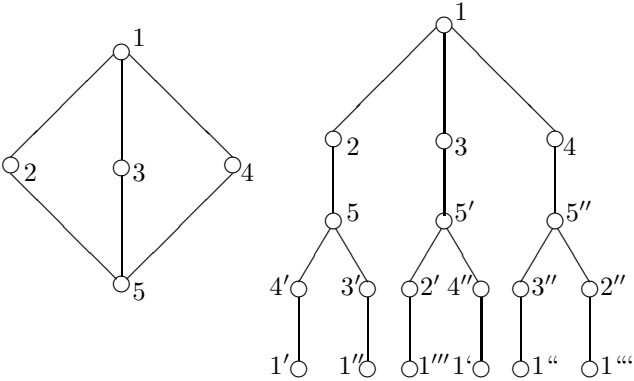

As in the works \cite{Roy1,Roy2,Weiss,Malioutov,Su,Su1}, the computation tree is an essential tool for convergence analysis of the message-passing algorithm on loopy graphs. We follow \cite{Roy1} for its construction. Given a loopy graph $G=(V,E)$ and a root node $r\in V$, its computation tree of depth $t>0$, denoted by $\mathbf{G}=(\mathbf{V}, \mathbf{E})$, is constructed iteratively. Denote the mapping from $\mathbf{V}$ to $V$ by $\sigma(\cdot)$. Without loss of generality, assume $r=1$. Placing node $r$ as the root node of $\mathbf{G}$, its child nodes correspond to all the neighbouring nodes in $N_r$. This forms the depth-1 computation tree. In each of the subsequent $(t-1)$ iterations, take each leaf node $i$ with parent node $j$. Then add all the neighbouring nodes in $N_{\sigma(i)} \backslash \sigma(j)$ in $G$ as the new child nodes. The edge set $\mathbf{E}$ is formed by connecting every child-parent pair, i.e., $(i,j)\in \mathbf{E}$ if and only if $(\sigma(i),\sigma(j))\in E$.  This construction is depicted in Fig.~\ref{fig:2-1} for $t=3$. 

The objective function $\mathbf{F}:\mathbb{R}^{|\mathbf{V}|} \rightarrow \mathbb{R}$ for $\mathbf{G}$, also in a pairwise separation form, is formed by taking the following:
\begin{itemize}
\item For each edge $(i,j)\in \mathbf{E}$, take $\mathbf{f}_{ij}(\cdot,\cdot):= f_{\sigma(i)\sigma(j)}(\cdot,\cdot)$ and  $\mathbf{J}_{i\rightarrow j}^{(0)}(\cdot) =J_{\sigma(i)\rightarrow \sigma(j)}^{(0)}(\cdot)$;
\item For each interior node (non-leaf node) $i\in \mathbf{V}$,  take $\mathbf{f}_i(\cdot):=f_{\sigma(i)}(\cdot)$; 
\item For each leaf node $i\in \mathbf{V}$ with parent node $j$,  take 
$$\mathbf{f}_i(\cdot):=f_{\sigma(i)}(\cdot)+\sum_{u\in N_{\sigma(i)}\backslash \sigma(j)}J_{u\rightarrow \sigma(i)}^{(0)}(\cdot).$$
\end{itemize}  
%Also define $\mathbf{J}_{i\rightarrow j}^{(0)}(\cdot) =J_{\sigma(i)\rightarrow \sigma(j)}^{(0)}(\cdot)$, etc. But whenever there is no confusion, we will take the convention that $x_i$ for $\mathbf{G}$ means $x_{\sigma(i)}$ for $G$. 

The {\em key} property of the computation tree is that, for the root node $r$ after $t$ iterations of message-passing, the estimate $x_r^{(t)}$ is identical for the original graph and the computation tree; see \cite{Roy1}. In addition, since $\mathbf{G}$ is a depth-$t$ tree graph, the optimal solution for minimising $\mathbf{F}(\cdot)$, denoted by $\tilde{x}$ in the sequel, is obtained by the message-passing algorithm after $t$ iterations; see \cite{Roy1}. Putting the two properties together, we have $x_r^{(t)} = \tilde{x}_r$.   

\subsection{Convergence Analysis of Message-Passing}

Consider the parameterised initial messages as in (\ref{eq:initial1}) with parameter $p\in[0,1]$ and use them to construct the depth-$t$ computation tree $\mathbf{G}$ with any root node $r\in V$. Denote the objective function for $\mathbf{G}$ by $\mathbf{F}(x,p)$ and the optimal solution by $\tilde{x}(p)$. In the sequel, we take the abuse of notation by naming $\sigma(i)$ as $i$ whenever there is no confusion. Similar to \cite{Roy1}, the first-order optimality conditions for $\mathbf{F}(x,p)$ are given as follows. For any interior node $i$ of $\mathbf{G}$, 
\begin{align}
&\nabla f_i(\tilde{x}_i(p))+\sum_{u\in N_i} \nabla_1 f_{iu}(\tilde{x}_i(p), \tilde{x}_u(p))=0. \label{eq:interior}
\end{align}
For any leaf node $i$ with parent node $j$, 
\begin{align}
&\nabla f_i(\tilde{x}_i(p))+\nabla_2 f_{ji}(\tilde{x}_j(p), \tilde{x}_i(p))\nonumber \\
+&\sum_{u\in N_i\backslash j} \nabla J_{u\rightarrow i}^{(0)}(\tilde{x}_i(p))+pC_{iu}=0. \label{eq:leaf}
\end{align}
Differentiating (\ref{eq:interior}) and (\ref{eq:leaf}) with respect to $p$ gives the following, similar to \cite{Roy1}.  For any interior node $i$ of $\mathbf{G}$, 
\begin{align}
& \left (\nabla^2 f_i(\tilde{x}_i)+\sum_{u\in N_i} \nabla_1^2 f_{iu}(\tilde{x}_i, \tilde{x}_u)\right )\nabla \tilde{x}_i \nonumber \\
&+ \sum_{u\in N_i} \nabla_{12} f_{iu}(\tilde{x}_i, \tilde{x}_u)\nabla \tilde{x}_u=0. \label{eq:interior1}
\end{align}
In the above, dependence of $\tilde{x}$ on $p$ is suppressed for convenience, and $\nabla \tilde{x}_i=\frac{d}{dp} \tilde{x}_i(p)$.  Similarly, for any leaf node $i$ with parent node $j$, 
\begin{align}
&\left (\nabla^2 f_i(\tilde{x}_i)+\nabla_2^2 f_{ji}(\tilde{x}_j, \tilde{x}_i)+\sum_{u\in N_i\backslash j} \nabla^2 J_{u\rightarrow i}^{(0)}(\tilde{x}_i)\right )\nabla \tilde{x}_i\nonumber \\
&+\nabla_{12} f_{ji}(\tilde{x}_j, \tilde{x}_i)\nabla \tilde{x}_j+\sum_{u\in N_i\backslash j}C_{iu}=0. \label{eq:leaf1}
\end{align}
Rewriting (\ref{eq:interior1}) and (\ref{eq:leaf1}) in a condensed form gives 
\begin{align}
\Gamma z &= h. \label{eq:Gzh}
\end{align}
In the above, the vector $z=\mathrm{col}\{\nabla \tilde{x}_i\}$. The matrix function  $\Gamma(\tilde{x})$ (or simply $\Gamma$) $=\{\gamma_{ij}\}$ is symmetric with entries specified below. If $i$ is an interior node of $\mathbf{G}$,   
\begin{align*}
\gamma_{ii} &= \nabla^2 f_i(\tilde{x}_i)+\sum_{u\in N_i} \nabla_1^2 f_{iu}(\tilde{x}_i, \tilde{x}_u)=\frac{\partial^2}{\partial x_i^2}F(\tilde{x}), \\
\gamma_{iu} &=\nabla_{12} f_{iu}(\tilde{x}_i, \tilde{x}_u), \ \ \forall u\in N_i.
\end{align*}
If $i$ is a leaf node of $\mathbf{G}$ with parent $j$, 
\begin{align*}
\gamma_{ii} &= \nabla^2 f_i(\tilde{x}_i)+\nabla_2^2 f_{ji}(\tilde{x}_j, \tilde{x}_i)+\sum_{u\in N_i\backslash j} \nabla^2 J_{u\rightarrow i}^{(0)}(\tilde{x}_i), \\
\gamma_{ij} &=\nabla_{12} f_{ji}(\tilde{x}_j, \tilde{x}_i).
\end{align*}
All other off-diagonal entries $\gamma_{iu}$ are zero. The vector $h=\mathrm{col}\{h_i\}$ is such that $h_i=0$ for every interior node $i$, and if $i$ is a leaf node with parent node $j$, 
\begin{align}
h_i&=-\sum_{u\in N_i\backslash j}C_{iu}.\label{eq:hi}
\end{align}
  
\begin{lem}\label{lem:3}
Under Assumptions~\ref{as:1}-\ref{as:2},  the matrix function $\Gamma(\tilde{x})$ is $(\lambda, \mathbf{w})$-scaled diagonally dominant for all $\tilde{x}\in \mathbb{R}^{|\mathbf{V}|}$, where $\mathbf{w}_i=w_{\sigma(i)}$ for all $i\in \mathbf{V}$.
\end{lem} 

\begin{IEEEproof}
The scaled diagonal dominance condition for an interior node $i$ of $\mathbf{G}$ is obvious by the construction of $\Gamma$ and Assumption~\ref{as:1}. For any leaf node $i$ of $\mathbf{G}$ with parent node $j$, 
\begin{align*}
&\lambda \mathbf{w}_i\gamma_{ii} - \mathbf{w}_j |\gamma_{ij}|\\
=&  \lambda w_i(\nabla^2 f_i(\tilde{x}_i)+\nabla_2^2 f_{ji}(\tilde{x}_j, \tilde{x}_i)+\sum_{u\in N_i\backslash j} \nabla^2 J_{u\rightarrow i}^{(0)}(\tilde{x}_i)) \\
&-w_j |\nabla_{12} f_{ji}(\tilde{x}_j, \tilde{x}_i)|\\
\ge& \lambda w_i(\nabla^2 f_i(\tilde{x}_i)+\nabla_2^2 f_{ji}(\tilde{x}_j, \tilde{x}_i)+\sum_{u\in N_i\backslash j}
\nabla_1^2f_{iu}(\tilde{x}_i,x_{u\rightarrow i}^{(0)}))\\
&-\sum_{u\in N_i\backslash j}w_u|\nabla_{12}f_{iu}(\tilde{x}_i,x_{u\rightarrow i}^{(0)})|-w_j|\nabla_{12}f_{ij}(\tilde{x}_i,\tilde{x}_j)|
\ge 0.
\end{align*}
The first inequality step above used Assumption~\ref{as:2}, and the last inequality step above used the scaled diagonal dominance property for $F(x)$. This completes the proof. 
\end{IEEEproof}

Now we are ready to prove Theorem~\ref{thm:2}. 

\begin{IEEEproof}
Take any node $\in V$ as the root node and construct the depth-$t$ computation tree $\mathbf{G}=(\mathbf{V},\mathbf{E})$ for any $t>0$. Define $\mathbf{W}=\mathrm{diag}({\mathbf{w}})$, $\Gamma_{\mathbf{w}}=\mathbf{W}^{-1}\Gamma\mathbf{W}$, $z_{\mathbf{w}}=\mathbf{W}^{-1}z$ and $h_{\mathbf{w}}=\mathbf{W}^{-1}h$. Then, (\ref{eq:Gzh}) can be rewritten as $\Gamma_{\mathbf{w}}z_{\mathbf{w}}=h_{\mathbf{w}}$. Further define $D=\mathrm{diag}(\Gamma)$ and $R_{\mathbf{w}}=I-D^{-1}\Gamma_{\mathbf{w}}$, we have $\Gamma_{\mathbf{w}}=D(I-R_{\mathbf{w}})$. By Lemma~\ref{lem:3}, it is clear that $\lambda^{-1}R_{\mathbf{w}}$ is diagonally dominant. It follows that $\|R_{\mathbf{w}}\|_{\infty}\le \lambda$, where the infinity norm for a matrix $A$ is defined to be
\begin{align*}
\|A\|_{\infty}:=\max_{i}\sum_{j} |a_{ij}|.
\end{align*}
Subsequently, the solution of $z_{\mathbf{w}}$ is given by 
\begin{align*}
z_{\mathbf{w}}&=(I-R_{\mathbf{w}})^{-1}D^{-1}h_{\mathbf{w}}=\sum_{s=0}^{\infty} R_{\mathbf{w}}^s b_{\mathbf{w}}  
\end{align*}
with $b_{\mathbf{w}}:=D^{-1}h_{\mathbf{w}}=D^{-1}\mathbf{W}^{-1}h$. In particular, the solution for the root node $r$ is given by 
\begin{align*}
z_r&=\mathbf{w}_r\sum_{s=0}^{\infty} e_r^TR_{\mathbf{w}}^s b_{\mathbf{w}},
\end{align*}
where $e_r$ is the column vector with 1 in entry $r$ and zero everywhere else. 

Following the {\em walk sum} argument in \cite{Malioutov}, $e_r^TR_{\mathbf{w}}^sb_{\mathbf{w}}$ is the walk sum of all the walks of length $s$ in $\mathbf{G}$ from any node to node $r$. As pointed out in \cite{Roy1}, since $h_i$ is nonzero only for the leaf nodes in the $t$-th depth of $\mathbf{G}$, the sum for $z_r$ above only needs to start from $s=t$. That is,  
\begin{align*}
\frac{z_r}{w_r}&=e_r^T\sum_{s=t}^{\infty} R_{\mathbf{w}}^s b_{\mathbf{w}}.
\end{align*}
It follows that 
\begin{align*}
\frac{|z_r|}{w_r}
&\le \|R_{\mathbf{w}}^t\|_{\infty} \sum_{s=0}^{\infty} \|R_{\mathbf{w}}^s\|_{\infty} \|b_{\mathbf{w}}\|_{\infty} \\
&\le \lambda^t\sum_{s=0}^{\infty} \lambda^s \max_{i\in \mathbf{V}_L} w_i^{-1}D_i^{-1}|h_i|  \\
&\le \frac{\lambda^t}{1-\lambda} \max_{i\in V} w_i^{-1}D_i^{-1}|h_i|. 
\end{align*}
In the above, $\mathbf{V}_L$ denotes the set of depth-$t$ leaf nodes in $\mathbf{G}$. 

Note that the above bound for $z_r(\rho)$ is valid for all $0\le \rho\le 1$.  Using the mean value theorem on $\tilde{x}(p)$, we get
\begin{align*}
\tilde{x}_r(1) - \tilde{x}_r(0) &= \nabla \tilde{x}_r(\breve{p}) (1-0) =z_r(\breve{p})
\end{align*}
with some intermediate value $\breve{p}\in [0, 1]$. Using $\tilde{x}_r(1)=x_r^{\star}$ (Lemma~\ref{lem:2}) and $\tilde{x}_r(0)=x_r^{(t)}$, it follows that 
\begin{align}
\frac{|x_r^{(t)} - x_r^{\star}|}{w_r} &\le \frac{\lambda^t}{1-\lambda} \max_{i\in V} w_i^{-1}D_i^{-1}|h_i|. \label{eq:zr}
\end{align}

Using (\ref{eq:hi}) and noting $D_i=\gamma_{ii}=\frac{\partial^2}{\partial x_i^2}F(\tilde{x})$, we have 
\begin{align}
D_i^{-1}|h_i|&\le \frac{\sum_{u\in N_i}|\nabla_1 f_{iu}(x_i^{\star},x_u^{\star})-\nabla J_{u\rightarrow i}^{(0)}(x_i^{\star})|}{\min_{x\in \mathbb{R}^n} \frac{\partial^2}{\partial x_i^2}F(x)}.\label{eq:Di1}
\end{align}
Hence, (\ref{eq:rate}) is verified.  

To verify (\ref{eq:rate1}), we first note that the special initial messages $J_{i\rightarrow j}^{(0)}(x_j)=f_{ji}(x_j,x_i^{(0)})$ also satisfy Assumption~\ref{as:2}. Using the mean value theorem, we get 
\begin{align*}
&\nabla_1 f_{iu}(x_i^{\star},x_u^{\star})-\nabla J_{u\rightarrow i}^{(0)}(x_i^{\star})\\
 =&\nabla_1 f_{iu}(x_i^{\star},x_u^{\star})-\nabla_1 f_{iu}(x_i^{\star}, x_u^{(0)}) \\
=& \nabla_{12} f_{iu}(x_i^{\star},\breve{x}_u)(x_u^{\star}-x_u^{(0)}) 
\end{align*}
for some intermediate value $\breve{x}_u$. Therefore, 
\begin{align*}
&\sum_{u\in N_i}|\nabla_1 f_{iu}(x_i^{\star},x_u^{\star})-\nabla J_{u\rightarrow i}^{(0)}(x_i^{\star})|\\
 \le &\sum_{u\in N_i}w_u|\nabla_{12} f_{iu}(x_i^{\star},\breve{x}_u)|\frac{|x_u^{\star}-x_u^{(0)}|}{w_u}\\
 \le & \sum_{u\in N_i}w_u|\nabla_{12} f_{iu}(x_i^{\star},\breve{x}_u)|\max_{v\in V}\frac{|x_v^{\star}-x_v^{(0)}|}{w_v} \\
 \le & \max_{x\in \mathbb{R}^n} \lambda w_i\frac{\partial^2}{\partial x_i^2}F(x)\max_{v\in V}\frac{|x_v^{\star}-x_v^{(0)}|}{w_v}. 
\end{align*}
The last step above used the scaled diagonal dominance property. Using the above in (\ref{eq:zr})-(\ref{eq:Di1}), we get 
\begin{align*}
\frac{|x_r^{(t)} - x_r^{\star}|}{w_r}&\le \frac{\lambda^{t+1}}{1-\lambda} \max_{i\in V} \frac{\max_{x\in \mathbb{R}^n}\hspace{-1mm} \frac{\partial^2}{\partial x_i^2}F(x)}{\min_{x\in \mathbb{R}^n}\hspace{-1mm}\frac{\partial^2}{\partial x_i^2}F(x)}\max_{v\in V}\hspace{-1mm}\frac{|x_v^{\star}-x_v^{(0)}|}{w_v}\\
&\le \frac{\lambda^{t+1}}{1-\lambda} M \max_{v\in V}\frac{|x_v^{\star}-x_v^{(0)}|}{w_v}.
\end{align*}
This completes the proof of (\ref{eq:rate1}).
\end{IEEEproof}

\section{Message-Passing for Quadratic Optimisation}\label{sec5}

In this section, we specialise the results in Section~\ref{sec3} to quadratic optimisation where the objective function becomes 
\begin{align}
F(x) &= \frac{1}{2}x^TAx -b^Tx, \label{eq:quadratic}
\end{align}
for some symmetric and positive-definite matrix $A=\{a_{ij}\}\in \mathbb{R}^{n\times n}$ and vector $b=[b_1, b_2, \ldots, b_n]^T\in \mathbb{R}^n$. It is obvious that such a $F(x)$ has a natural pairwise separation (\ref{eq:ps}) with 
\begin{align}
f_i(x_i) &= \frac{1}{2}a_{ii} x_i^2-b_ix_i, \ \ f_{ij}(x_i,x_j)= a_{ij}x_ix_j,  \label{eq:q1}
\end{align}
but this is not a pairwise convex separation (unless all $a_{ij}=0$). Although it is known \cite{Malioutov,Roy1} that $F(\cdot)$ being pairwise convex separable if and only if $A$ is scaled diagonally dominant, finding a corresponding scaling vector $w$ and re-parameterising $F(x)$ in a pairwise convex separation would be an optimisation task on its own. 
%To our knowledge, no distributed algorithm is seen in the literature for this task. 

By specialising Theorems~\ref{thm:1} and \ref{thm:2} to the quadratic case, we provide a very simple convergence property of the message passing algorithm without requiring a known pairwise convex separation form. That is, we can work directly using the pairwise separation form (\ref{eq:q1}). This makes the message-passing algorithm in line with \cite{Weiss,Malioutov,Roy2}, but the new contribution here is that a convergence rate is explicitly presented. 

\begin{assumption}\label{as:3}
\begin{itemize}
\item  The matrix $A$ is $(\lambda-w)$-scaled diagonally dominant for some $0<\lambda<1$ and positive $w\in \mathbb{R}^n$. 
\item The initial messages are chosen as,  $\forall \ i\in V, j\in N_i$, 
\begin{align}
J_{i\rightarrow j}^{(0)}(x_j) &= \frac{1}{2}\alpha_{i\rightarrow j}^{(0)} x_j^2-\beta_{i\rightarrow j}^{(0)}x_j, \label{eq:Jij3}
\end{align}
with any constants $\alpha_{i\rightarrow j}^{(0)}$ and $\beta_{i\rightarrow j}^{(0)}$ satisfying 
\begin{align}
\alpha_{i\rightarrow j}^{(0)}\ge -\rho\frac{w_i}{w_j}|a_{ji}| \label{eq:alphaij}
\end{align}
for some $0\le \rho<\lambda^{-1}$.
\end{itemize}
\end{assumption}
Note that $\beta_{i\rightarrow j}^{(0)}$ has no constraints and that  the initial estimates are actually ``hidden'' in $\beta_{i\rightarrow j}^{(0)}$.

\subsection{Well-posedness for Quadratic Optimisation}
  
Specialising Theorem~\ref{thm:1} to the quadratic case, we get the following result:
\begin{thm}\label{thm:3}
Under Assumption~\ref{as:3}, the following properties hold for all $i\in V, j\in N_i$ and $t\ge0$:
\begin{itemize}
\item[P1':] The functions $\gamma_{i}^{(t)}(\cdot)$ are strictly convex and the functions $g_{ij}^{(t)}(\cdot, \cdot)$ are strictly convex with respect to the first variable, hence the message-passing algorithm is well posed;
\item[P2':] The new messages are given by
\begin{align}
J_{i\rightarrow j}^{(t+1)}(x_j)&= \frac{1}{2}\alpha_{i\rightarrow j}^{(t+1)}x_j^2-\beta_{i\rightarrow j}^{(t+1)}x_j, \label{eq:Jij4}
\end{align}
where
\begin{align}
\alpha_{i\rightarrow j}^{(t+1)}&= -\frac{(a_{ji})^2}{a_{i \rightarrow j}^{(t)}}, \ \ \beta_{i\rightarrow j}^{(t+1)}=-\frac{a_{ji}b_{i\rightarrow j}^{(t)}}{a_{i \rightarrow j}^{(t)}}  \label{eq:alphabeta}
\end{align}
with 
\begin{align}
a_{i\rightarrow j}^{(t)} &= a_{ii} +\hspace{-3mm}\sum_{u\in N_i\backslash j} \alpha_{u\rightarrow i}^{(t)}, \ 
b_{i\rightarrow j}^{(t)} =b_i + \hspace{-3mm} \sum_{u\in N_i\backslash j} \beta_{u\rightarrow i}^{(t)}.   \label{eq:abij}
\end{align}
\item[P3':] $\lambda w_ia_{i\rightarrow j}^{(t+1)}\ge w_j|a_{ji}|>0$ and $0>\alpha_{i\rightarrow j}^{(t+1)}\ge - \frac{w_i}{w_j}|a_{ji}|$. 
\end{itemize}
\end{thm}
  
\begin{rem}\label{rem:x}
The property (P3') in Theorem~\ref{thm:3} means that all the new messages $J_{i\rightarrow j}^{(t+1)}$ are actually {\em concave}, even if the initial messages $J_{i\rightarrow j}^{(0)}$ may be convex.  In contrast, the message-passing algorithm in \cite{Roy1} requires that all the messages need to remain convex. The reason for the new messages to be concave in our case is due to the fact that the given pairwise separation (\ref{eq:q1}) has non-convex $f_{ij}(\cdot,\cdot)$, whereas \cite{Roy1} requires the given $f_{ij}(\cdot,\cdot)$ to be convex. 
\end{rem}

\subsection{Convergence Properties}

\begin{thm}\label{thm:4}
Under Assumption~\ref{as:3}, the following holds for the message-passing algorithm when applied to (\ref{eq:quadratic}): 
\begin{align}
\frac{|x_r^{(t)}-x_r^{\star}|}{w_r} &\le \frac{\lambda^{t}}{1-\lambda}\max_{i\in V} \frac{\sum_{u\in N_i}|a_{iu}x_u^{\star}-\alpha_{i\rightarrow j}^{(0)}x_i^{\star}+\beta_{u\rightarrow i}^{(0)}|}{w_ia_{ii}},\label{eq:rateq}
\end{align}
 for every $r\in V$ and $t>0$. If we choose $J_{i\rightarrow j}^{(0)}(x_j)=a_{ji}x_jx_i^{(0)}$ for some initial $x_i^{(0)}$ (i.e., $\alpha_{i\rightarrow j}^{(0)}=0, \beta_{i\rightarrow j}^{(0)} = -a_{ji}x_i^{(0)}$) for all $i\in V, j\in N_i$,  then
\begin{align}
\frac{|x_r^{(t)}-x_r^{\star}|}{w_r} &\le \lambda^{t+1}\frac{1}{1-\lambda}\max_{v\in V}\frac{|x_v^{(0)}-x^{\star}|}{w_v}.\label{eq:rateq1}
\end{align}
  
\end{thm}

\begin{IEEEproof}
Note that, in the quadratic case, $M$ as defined in (\ref{eq:M}) equals 1. The proof follows directly from the properties in Theorem~\ref{thm:3} and the convergence results in Theorem~\ref{thm:2}. 
\end{IEEEproof}

\begin{rem}\label{rem:q}
The work of \cite{Roy2} on message-passing for quadratic optimisation showed that the message-passing algorithm converges asymptotically under the assumption that the initial messages conform to a {\em convex-dominated decomposition}. In comparison, Theorem~\ref{thm:4} gives an explicit convergence rate bound.
\end{rem}

\section{Conclusions}\label{sec6}

In this paper, we followed the technique in \cite{Roy1} and improved their convergence analysis of message-passing for convex optimisation by removing the requirement that the given pairwise components need to be convex. We have also provided a tighter and simpler bound for the convergence rate. Our convergence rate bound for the quadratic optimisation case is also new. We have not considered the so-called {\em asynchronous} version of the algorithm, but similar results are expected to hold, as done in \cite{Weiss,Roy1,Roy2,Malioutov}. We stress that the bounds on the convergence rate are still conservative in general. Further work is needed to find even better bounds. How to relax the diagonal dominance assumption will be of future interest as well.

\end{document}